\documentclass[12pt]{article}
\usepackage{amssymb, amstext, amsmath, amsthm, amscd, graphicx, amsfonts}
\usepackage{latexsym}
\newtheorem{teorema}{Theorem}
\newtheorem{lema}[teorema]{Lemma}
\newtheorem{remark}[teorema]{Remark}
\newtheorem{observation}[teorema]{Observation}

\begin{document}

\global\newcount\numsec\global\newcount\numfor

\def\sqr#1#2{{\vcenter{\vbox{\hrule height.#2pt
\hbox{\vrule width.#2pt height#1pt \kern#1pt \vrule
width.#2pt}\hrule height.#2pt}}}}
\def\square{\mathchoice\sqr34\sqr34\sqr{2.1}3\sqr{1.5}3}

\let\ni\noindent
\let\ciao=\bye \def\fiat{{}}
\def\pagina{{\vfill\eject}} \def\\{\noindent}
\def\bra#1{{\langle#1|}} \def\ket#1{{|#1\rangle}}
\def\media#1{{\langle#1\rangle}} \def\ie{\hbox{\it i.e.\ }}
\let\ii=\int \let\ig=\int \let\io=\infty \let\i=\infty

\let\dpr=\partial \def\V#1{\vec#1} \def\Dp{\V\dpr}
\def\oo{{\V\o}} \def\OO{{\V\O}} \def\uu{{\V\y}} \def\xxi{{\V \xi}}
\def\xx{{\V x}} \def\yy{{\V y}} \def\kk{{\V k}} \def\zz{{\V z}}
\def\rr{{\V r}} \def\zz{{\V z}} \def\ww{{\V w}}
\def\Fi{{\V \phi}}
\def\an{{\cal A}}
\let\Rar=\Rightarrow
\let\rar=\rightarrow
\let\LRar=\Longrightarrow

\def\lh{\hat\l} \def\vh{\hat v}

\def\ul#1{\underline#1}
\def\ol#1{\overline#1}

\def\ps#1#2{\psi^{#1}_{#2}} \def\pst#1#2{\tilde\psi^{#1}_{#2}}
\def\pb{\bar\psi} \def\pt{\tilde\psi}

\def\E#1{{\cal E}_{(#1)}} \def\ET#1{{\cal E}^T_{(#1)}}
\def\RR{{\cal R}}\def\SS{{\cal S}} \def\HH{{\cal
H}}\def\GG{{\cal G}}
\def\PP{{\cal P}} \def\AA{{\cal A}}
\def\BB{{\cal B}}\def\FF{{\cal F}}
\def\GI{{\mathbb{G}}}\def\V{{\mathbb{V}}}\def\E{{\mathbb{E}}}
\def\tende#1{\vtop{\ialign{##\crcr\rightarrowfill\crcr
              \noalign{\kern-1pt\nointerlineskip}
              \hskip3.pt${\scriptstyle #1}$\hskip3.pt\crcr}}}
\def\otto{{\kern-1.truept\leftarrow\kern-5.truept\to\kern-1.truept}}
\def\arm{{}}
\font\bigfnt=cmbx10 scaled\magstep1

\def\PP{\mathbb{P}}\def\Pd{\mathbb{P}}
\def\LL{\mathbb{L}}\def\Ld{\mathbb{L}}
\def\NN{\mathbb{N}}\def\Nd{\mathbb{N}}
\def\EE{\mathbb{E}}\def\Ed{\mathbb{E}}
\def\VV{\mathbb{V}}\def\Vd{\mathbb{V}}
\def\ZZ{\mathbb{Z}}\def\Lad{{\mathbb{Z}_+^2}}
\def\Rd{{\mathbb{R}^2}}
\def\0{\emptyset}
\def\F{{\cal F}}
\def\v{\vskip.1cm}
\def\vv{\vskip.2cm}
\def\vvv{\vskip.3cm}
\def\vvvv{\vskip.4cm}

\date{}

\title{Linear Lower Bounds for $\delta_c(p)$ for a Class of $2D$ Self-Destructive Percolation Models}

\author{J. van den Berg \footnote {Centrum voor Wiskunde en Informatica, Amsterdam,
the Netherlands; J.van.den.Berg@cwi.nl; part of vdB's research is funded by the Dutch
BSIK/BRICKS project.} and B.N.B. de Lima \footnote{ Universidade
Federal de Minas Gerais, Belo Horizonte, Brazil; bnblima@mat.ufmg.br; partially
supported by CNPq.}} \maketitle
\begin{abstract}
The self-destructive percolation model is defined as follows:
Consider percolation with parameter $p > p_c$. Remove the infinite
occupied cluster. Finally, give each vertex (or, for bond
percolation, each edge) that at this stage is vacant, an extra
chance $\delta$ to become occupied. Let $\delta_c(p)$ be the
minimal value of $\delta$, needed to obtain an infinite occupied
cluster in the final configuration. This model was introduced some
years ago  by van den Berg and Brouwer. They showed that, for the
site model on the square lattice (and a few other $2D$ lattices
satisfying a special technical condition) that
$\delta_c(p)\geq\frac{(p-p_c)}{p}$. In particular, $\delta_c(p)$
is at least linear in $p-p_c$.

Although the arguments used by van den Berg and Brouwer look quite rigid, we show that they
can be suitably modified to obtain
similar linear lower bounds for $\delta_c(p)$ (with $p$ near $p_c$)
for a much larger class of
$2D$ lattices, including bond
percolation on the square and triangular lattices, and site percolation on the star lattice (or matching
lattice) of the square lattice.
\end{abstract}

\vv
{\footnotesize

\v
}
\vvvv
\section{Introduction}

Some years ago van den Berg and Brouwer, motivated by the study of forest-fire processes,
introduced the self-destructive percolation model (see \cite{BB} and
\cite{B}).
This model can be described in a few steps as follows.

Let $G=(\VV,\EE)$ be a connected, infinite, locally finite graph. The first step is to
perform independent site percolation with parameter $p$ on this
graph: we declare each vertex occupied with probability $p$
and vacant with probability $1-p$, independent of the other vertices.
If $U, W \subset \VV$,
we use the notation $\{U\leftrightarrow W\}$ for the
event that there is a path of occupied vertices connecting some
vertex of $U$ to some vertex of $W$,  and $\{U\leftrightarrow
\infty\}$ for the event that there exists an infinite path of
occupied vertices starting from some vertex in $U$. If $G=(\VV,\EE)$ is
transitive (i.e., without lost of generality we can choose
any vertex to be the origin), let $\theta (p)$ be the probability
that the origin belongs to an infinite occupied cluster.

Since $\theta(p)$ is monotone in $p$, there is a critical value $p_c\in [0,1]$ \\ such
that $\theta (p)>0$ if $p\in (p_c,1]$ and $\theta (p)=0$ if $p\in
[0,p_c)$. It is well-known (see \cite{G} for this and other background results on
percolation) that $\theta (p_c)=0$ for a large class
of graphs including the $d$-dimensional hypercubic lattice, for $d=2$ or
$d\geq 19$, the triangular and hexagonal lattices.

Now suppose that, by some `catastrophe', the infinite occupied
clusters are destroyed, that is, each vertex that belongs to an
infinite occupied cluster becomes vacant.
After the
catastrophe we give each vacant vertex an extra chance to become
occupied; more precisely, each vacant vertex becomes occupied
with probability $\delta$, independent of everything else. Let
$P_{p,\delta}$ be the distribution of the final configuration and
$\theta (p,\delta):=P_{p,\delta}(0\leftrightarrow\infty)$.

An equivalent but more formal (and often more convenient) description of the
self-destructive percolation model is the
following: let $X_i,\ i\in\VV$ be a sequence of i.i.d. 0-1 random
variables with parameter $p$, and let $Y_i,\ i\in\VV$ be another
sequence of i.i.d. 0-1 random variables with parameter $\delta$.
(Here, we can interpret 0 as vacant and 1 as occupied).
Moreover we take the sequence of $Y_i$'s independent of the
sequence of $X_j$'s. Let $X_i^*,\ i\in\VV$ be defined by
\begin{equation}
X^*_i=\left\{
\begin{array}
[c]{l}%
1,\mbox{  if } X_i=1 \mbox{ and there is no infinite $X$-occupied path from $i$} \\
0,\mbox{ otherwise}
\end{array}\right.
\end{equation}

Here, by $X$-occupied path we mean a path on which each vertex $j$ has
$X_j=1$. Finally, we define $Z_i:=X_i^*\vee Y_i$. This sequence
$Z_i,\ i\in\VV$ is the final configuration and the measure
$P_{p,\delta}$ is its distribution. Analogously (with obvious modifications of the definitions) one can
define the self-destructive bond percolation model.

Monotonicity in $\delta$ obviously implies that
for each $p\in [p_c,1]$ there exists a
$\delta_c(p)\in [0,1]$ such that $\theta (p,\delta )=0$ for all
$\delta\in [0,\delta_c(p))$ and $\theta (p,\delta )>0$ for all
$\delta\in (\delta_c(p),1]$.
It is also clear that if $\delta > p_c$, then $\theta
(p,\delta)>0$ for all $p\in [0,1]$; that is,
\begin{equation} \label{triv-bound}
\forall p \,\, \delta_c(p) \leq p_c.
\end{equation}

As to lower bounds for $p_c$, Proposition 3.1 of
\cite{BB} and Proposition 2.3.2 in \cite{B} show that for site
percolation on the square (and the triangular and the honeycomb) lattice \\
$\delta_c(p)$ is at least linear in $p-p_c$. More precisely it was shown for
these models that

\begin{equation}\label{bound}
\mbox{If } p(1-\delta)>p_c, \mbox{ then }  \theta (p,\delta)=0,
\end{equation}
which is obviously equivalent to

\begin{equation}\label{bound2}
\delta_c(p)\geq\frac{p-p_c}{p}, \,\, p \geq p_c. \\
\end{equation}

(In fact, it
is conjectured in that paper that for self-destructive percolation on these lattices $\delta_c(p)$ does not
go to $0$ as $p \downarrow p_c$, but a proof or disproof of that conjecture seems out of reach at the moment).

Using totaly different arguments it has been proved that
\eqref{bound2} also holds for site percolation on the binary tree
(see Theorem 2.6.1 in \cite{B} or Theorem 5.1 in \cite{BB}). This raises the question
how general \eqref{bound2}, or, at least, the following weaker property, are:

\begin{equation}\label{weakest}
\exists \hat p > p_c \,\, \exists C > 0 \mbox{ s.t. } \,\, \forall p\in[p_c,\hat p] \,\,
\delta_c(p) \geq C (p-p_c).
\end{equation}


As mentioned in Remark (ii) just after the above mentioned Proposition 3.1 in
\cite{BB}, the argument of this Proposition does not work
for bond percolation on the square lattice.
However, in Section 2 we will show that the statement of that Proposition is true for that model. We do this by
refining the argument in \cite{BB}.

On the other hand, it is easy to see that \eqref{bound2} can not hold for the bond model on the
triangular lattice, the site model on the matching (or star)
lattice of the square lattice, and more generally any percolation model with $p_c$
smaller than $1/2$: For such models $(p - p_c)/ p$ is clearly larger than
$1/2$ (and hence larger than $p_c$) for $p$ sufficiently close to 1.
For such $p$ \eqref{bound2} would contradict \eqref{triv-bound}.

Nevertheless, still further refinement of the arguments shows that for a
large class of 2D percolation models,
including the above mentioned bond model on the triangular lattice and the site model on the
star lattice of the square lattice, \eqref{weakest} does hold.
This is done in Section 3.

This paper focuses mainly on transitive graphs. However, related to Remark \ref{rem-gen} at the end
of Section 3, we make the
following comment on non-transitive graphs. On such graphs $P_p(v \leftrightarrow \infty)$ may depend on
the vertex $v$, and for this reason the notation $\theta(p;v)$ is used. However (as is well-known), by positive
association (FKG)
$\theta(p;w) \geq P_p(v \leftrightarrow w) \theta(p;v)$ for all vertices $v, w$.
Therefore $\theta(p;v)$ is positive if and only if
$\theta(p;w)$ is positive, and hence $p_c$ does not depend on $v$. For similar reasons (see \cite{BB} for
positive association for the self-destructive percolation model) $\delta_c(p)$ does not
depend on $v$ either.
\vvvv
\section{Bond percolation on the Square Lattice}

\vvv

Let $\LL^2=(\ZZ^2,\EE (\LL^2))$ be the square lattice, where $$\EE
(\LL^2)=\{\langle x,y\rangle ; \|x-y\|_1=1\}$$ and let
$G_{cb}=(\ZZ^2,\EE_{cb})$ be the chess-board lattice, where
$$\EE_{cb}=\{\langle x,y\rangle ;
\|x-y\|_1=1\}$$
$$\cup\ \{\langle (x_1,x_2),(y_1,y_2)\rangle
;(y_1,y_2)=(1,1)+(x_1,x_2)\mbox{ and }x_1+x_2\mbox{ is even }\}$$
$$\cup\ \{\langle (x_1,x_2),(y_1,y_2)\rangle
;(y_1,y_2)=(1,-1)+(x_1,x_2)\mbox{ and }x_1+x_2\mbox{ is odd }\}$$

It is well-known that bond percolation on the square lattice is
equivalent to site percolation on its covering graph (see section
2.5 in \cite{K} or section 1.6 in \cite{G}), the chess-board
lattice.

The next result is an extension of Proposition 3.1 in \cite{BB}
for bond percolation on the square lattice.
For the proof of that proposition it was essential that
the lattice under consideration
is a subgraph of its matching lattice. (For the notion `matching lattice', see e.g. Section 3.1
in \cite{G}. In the particular case of the square lattice, the matching lattice is obtained from the
square lattice by adding, in each face, the two diagonals as extra edges; the matching lattice of the
triangular lattice is the tringular lattice itself: it is self-matching).
The chess-board lattice does not have this property. This is why
the proof of \cite{BB} does not work for the site model on that lattice (and hence for bond percolation
on the square
lattice). However, the chess-board lattice is a translation (`along an edge')
of its matching lattice, and we will exploit this property to modify the proof
of the above mentioned Proposition 3.1 in \cite{BB}.

\begin{teorema} \label{square} For the self-destructive site
percolation model on the chess-board lattice (or, equivalently, bond
percolation on the square lattice), it holds that if $p(1-\delta
)>p_c$ then $\theta(p,\delta )=0$. Hence,
$\delta_c(p)\geq\frac{p-p_c}{p}$.
\end{teorema}

\begin{proof}
Given $v\in\ZZ^2$, we will use the notation $\tilde{v}:=v+(1,0)$.
Recall the sequences of 0-1 valued random variables, $X_v,Y_v,Z_v \
v\in\ZZ^2$, introduced in section 1. We color each vertex
$v\in\ZZ^2$ red if $X_v=1$ and $Y_{\tilde{v}}=0$. Then, each vertex
$v$ will be red with probability $p(1-\delta )$, independently of
the other vertices. Since, $p(1-\delta )>p_c$ it follows from
ordinary site percolation (see section 11.8 in \cite{G}), that, a.s.,  there is
an infinite red cluster, and this cluster contains a circuit around the
origin.

Let $\gamma$ be such a
red circuit.
Define $$\tilde{\gamma}:=\gamma + (1,0).$$
Note that $\tilde{\gamma}$ is a circuit in the matching graph.
Let $\tilde{v}\in\tilde{\gamma}$.
By construction,
$Y_{\tilde{v}}=0$. Moreover, since $\tilde{v}$ is a neighbor of the infinite
$X$-occupied cluster, $X^*_{\tilde{v}} = 0$. Hence $Z_{\tilde{v}}=0$. Summarizing we have that,
almost surely, there is a
$Z$-vacant circuit in the matching graph which surrounds or contains the
origin. Hence $\theta(p,\delta )=0$.
\end{proof}

\vvv

\vvvv
\section{Other  $2D$ lattices}

Note that in the proof of Theorem \ref{square} (as in that of
Proposition 3.1 of \cite{BB}) the definition of {\em red} vertices
(or edges) was done in such a way that each vertex (edge) is red
independently of the other vertices (edges). In the current
Section the color of a vertex will involve the $Y$ values of all
its neighbours. This strategy, which will be used to show that
\eqref{weakest} is true for a large class of $2D$ lattices, leads
to dependencies which somewhat complicate the analysis.

%

The main result of this section, Theorem \ref{general} below, is stated for three well-known lattices,
but, as we shall point out in Remark \ref{rem-gen}, the result (with practically the same proof) holds for a
large class of $2 D$ lattices.

\begin{teorema} \label{general}
For self-destructive site percolation on the matching lattice of the square or honeycomb lattice, or self-destructive
bond percolation on the triangular lattice, the following holds:

There are $\hat p >p_c$ and $C > 0$ such that
$$\forall p \in [p_c, \hat p] \, \delta_c(p) \geq C (p - p_c).$$
\end{teorema}

\begin{proof}
From now on $G = (\VV, \EE)$ denotes the matching lattice of the square or honeycomb lattice, or the
covering lattice of the triangular lattice, and $p_c$ the
critical probability for site percolation on $G$. (See, however, Remark \ref{rem-gen} below).

Recall that for the proofs in Section 2 we introduced a certain
colouring of the vertices, and that the colours were i.i.d. so
that we could compare the result of the coulouring with ordinary
percolation. In the current situation we will again define a
colouring, but now the colours are not independent. Nevertheless
it turns out that we again obtain a suitable comparison with
ordinary percolation. First some notation.

For each vertex $v\in\VV$, let $D_v:=\{u\in\VV ; \langle v,u
\rangle\in\EE\}$. Let $d : =|D_v|$ (for example, $d= 8$ for the
matching lattice of the square lattice).

Recall the sequences $X_v,Y_v,\ v\in\VV$ defined in Section 1.
Now define the sequence of 0-1 random variables $R_v,\
v\in\VV$ as
\begin{equation}
R_v=\left\{
\begin{array}
[c]{l}%
1,\mbox{  if } X_v=1 \mbox{ and } Y_u=0,\ \forall u\in D_v\\
0,\mbox{ otherwise}
\end{array}\right.
\end{equation}
If $R_v = 1$ we say that $v$ is $R-$occupied (or, simply, that $v$ is red).

Before we proceed with the proof of the theorem, we first state Observation \ref{obs} and
state and prove Lemma \ref{constant} below,
which  will be used later.

\medskip\noindent
\begin{observation} \label{obs}
Let $\gamma$ be a circuit in $G$. Then every path that starts in the interior of $\gamma$ and ends in the
exterior of $\gamma$ contains a vertex which has a neighbour on $\gamma$.
\end{observation}

\begin{lema}\label{constant}Let $\epsilon >0$. There is a constant $c_{\epsilon}$ such that
for the self-destructive site
percolation model on $G$ with parameters $0 < \delta \leq p_c$ and $p \in(p_c,  1 - \epsilon)$ the
following holds:\\
For every $v\in\VV$, every finite
subset of vertices $F\subset \VV$ and every colouring $(r_u, \,  u\in F)$ of $F$,
\[\frac{P(Y_v=1; R_u=r_u,\,  u\in F)}{P(Y_v=0; R_u=r_u,\, u\in
F)}\leq c_{\epsilon} \delta .\]Hence, $P(Y_v=1 |R_u=r_u,\, u\in F)\leq c_{\epsilon} \delta $.
\end{lema}

{\em Proof of Lemma} \ref{constant}: If $D_v\cap F =\emptyset$, the random variable $Y_v$
is independent of the sequence $R_u,\, u\in F$, and hence
\[\frac{P(Y_v=1; \, R_u=r_u, u\in F)}{P(Y_v=0; \, R_u=r_u, u\in F)}= \frac{\delta}{1-\delta} .\]

For the case $ F^\prime := D_v\cap F \neq\emptyset$, we consider two
subcases: If there is at least one $u\in F^\prime$ with $r_u=1$
then obviously $$P(Y_v=1, R_u=r_u,\ \forall u\in F)=0,$$
and it is easy to see that $P(Y_v =0; \, R_u = r_u, u \in F) > 0$.

If $r_u=0,\ \forall u\in F^\prime$ we have that
\[P(Y_v=1; \, R_u=r_u, u\in F)=P(Y_v=1)P(R_u=r_u, u\in
F \, | \, Y_v=1)\]
\[\leq\ \delta P(R_u=r_u, u\in F \setminus F^\prime \, | \, Y_v = 1 ) \]
\[ =\ \delta
P(R_u=r_u, u\in F \setminus F^\prime),\]
and
\[P(Y_v=0; \,
R_u=r_u, u\in F)=P(Y_v=0)P(R_u=r_u, u\in F \, | \, Y_v=0)\]
\[\geq\ (1-\delta)P(R_u=r_u, u\in F; \, X_z=0, z\in D_v \, | \, Y_v=0)\]
\[=\  (1-\delta)(1-p)^d P(R_u=r_u, u\in F\, | \, Y_v=0; \, X_z=0, z\in D_v)\]
\[=\ (1-\delta)(1-p)^d P(R_u=r_u, u\in F \setminus F^\prime).\]

Combining these two inequalities we have that
\[\frac{P(Y_v=1; \, R_u=r_u,  u\in F)}{P(Y_v=0; \, R_u=r_u, u\in
F)}\leq \frac{\delta}{(1-\delta)(1-p)^d}.\] So the claim of the lemma holds,
with the constant
$c_{\epsilon} = \frac{1}{(1-p_c) \, \epsilon^d}$. \\
{\em This completes the proof of Lemma} \ref{constant}.


Now we continue the proof of Theorem \ref{general}.
Let $\epsilon > 0$.
Suppose $p \in (p_c, 1 -  \epsilon)$, and let $\delta$ be such that $p(1-d c_{\epsilon} \delta) > p_c$,
where $c_{\epsilon}$ is the constant given in Lemma \ref{constant}.


Let $v\in\VV$, $F$ a finite set of vertices not containing $v$ and
$r_u\in\{0,1\},\ \forall u\in F$. We have
\[P(R_v=1 \, |\, R_u=r_u, u\in F)=P(X_v=1; \, Y_z=0, z\in D_v \, | \, R_u=r_u, u\in F)\]
\[=\ P(X_v=1)P(Y_z=0,\,  z\in D_v \, | \,  R_u=r_u, u\in
F)\]
\[\geq\ p(1-d c_{\epsilon} \delta),\] where in the inequality we
used Lemma \ref{constant} (and in the equality the fact that $X_v$ is independent of the
collection of random variables $\{ R_u, u \in F; Y_z, z \in D_v\}$).
As $p(1-d c_{\epsilon} \delta) > p_c$, the process $(R_v, v \in \VV)$ dominates an i.i.d. process
with parameter larger than $p_c$. Comparison with ordinary percolation shows that a.s. there is an infinite
$R$-occupied cluster which contains
a circuit around the origin.
Let $\gamma$ be such a circuit. Observe that, by the definition of the colourings, $\gamma$ belongs to an
infinite $X$-occupied cluster.
Define
$$\Gamma:=\cup_{v\in\gamma }D_v.$$
For each $w \in \Gamma$ we have either $X_w = 1$, in which case (by the above observation) $w$ belongs
to an infinite $X-$open cluster, or we have $X_w = 0$. In both cases $X_w^* = 0$. Since also $Y_w = 0$ for all
$w \in \Gamma$, we have $Z_w = 0$ for all $w \in \Gamma$. By Observation \ref{obs} we now conclude that
there is no infinite  $Z-$open path starting in $O$.
So we have proved that
for all $\epsilon > 0$ and all $p < 1-\epsilon$, it holds that
$\theta(p,\delta )=0$ if $p(1-d \,  c_{\epsilon} \, \delta) > p_c$; that is,

\begin{equation} \label{stronger}
 \delta_c(p) \geq \frac{p - p_c}{p \, d \, c_{\epsilon}}, \,\, p < 1 - \epsilon.
\end{equation}
This completes the proof of Theorem \ref{general}.
\end{proof}

\begin{remark} \label{rem-stronger}
Note that in fact we have proved something stronger than the claim in the theorem,
namely that for every $\epsilon > 0$ there is a $c_{\epsilon} > 0$ such that \eqref{stronger}
holds for all $p < 1 - \epsilon$.

This can be extended to the result that there is a $C > 0$ such that $\delta_c(p) \geq C (p-p_c$ for
all $p > p_c$. Since the $c_{\epsilon}$ in \eqref{stronger} goes to $\infty$ as $\epsilon$ goes to $0$, this result
does not follow immediately.
However, to get the extension,
it suffices to show that there is an $\epsilon >0$ such that $\delta_c(p)$ is bounded
away from $0$ for $p > 1-\epsilon$. Or, equivalently, that there are $\hat p < 1$ and $\delta >0$ such that

\begin{equation} \label{finclaim}
\theta(p,\delta) = 0 \mbox{ for all } p > \hat p.
\end{equation}
This can be proved by (e.g.) arguments very similar to those used in the proof that (ii) implies (i) in
Theorem 5.1 in \cite{BBV}.
Artem Sapozhnikov (private communication) has pointed out to us that further refinements of such arguments show that, on the $d$-dimensional cubic lattice, $\delta_c(p) \rightarrow p_c$
as $p \rightarrow 1$.
\end{remark}

%
%
%
\smallskip\noindent
\begin{remark} \label{rem-gen}
Theorem \ref{general} holds for a large class of 2D lattices. Essentially we only used that for supercritical
percolation the infinite cluster a.s. contains a circuit around $O$, which satisfies Observation \ref{obs}.
This property holds for site
percolation on the lattices belonging to the
family in Theorem 12.1 in Kesten's book \cite{K}. Informally speaking, this
family consists of lattices which belong to a pair of matching lattices
with certain periodicity and reflection symmetry properties (but which
are not necessarily transitive).
\end{remark}

\vvvv

\noindent {\bf Acknowledgments.} This work was done during de
Lima's sabbatical leave at CWI. He would like to thank CWI for the
hospitality and CAPES (Brazi\-li\-an Ministry of Education) for
the support during this period.


\begin{thebibliography}{99}

\bibitem{BB} van den Berg, J. and Brouwer, R., Self-Destructive Percolation,
{\em Random Structures and Algorithims} {\bf 24}, 480-501 (2004).

\bibitem{BBV} van den Berg, J., Brouwer, R. and V\'agv\"olgyi, B., {\em Box-crossings and
continuity results for self-destructive percolation in the plane}, to appear in Proceedings of XEBP
(Tenth Brazilian School of Probability, eds. V. Sidoravicius and M.E. Vares),  Birkh\"auser.

\bibitem{B}  Brouwer, R., {\em Percolation, forest-fires and monomer-dimers.} PhD Thesis,
Amsterdam, 2005.

\bibitem{G} Grimmett G., \emph{Percolation}, 2nd edition, Springer-Verlag, Berlin, 1999.


\bibitem{K} Kesten, H., \emph{Percolation Theory for
Mathematicians}, Birkh\"auser, Boston, 1982.

\bibitem{S} Sapozhnikov, A., Private Communication.

\end{thebibliography}
\end{document}